\documentstyle[12pt]{article}
\begin{document}

\newcommand{\R}{{\bf R}}
\newcommand{\N}{{\bf N}}
\newcommand{\Z}{{\bf Z}}
\begin{center}
Exponentiation and Euler measure \\
$ $ \\
James Propp \\
University of Wisconsin \\
April 29, 2001 \\
Revised May 29, 2002 \\
\end{center}

\begin{center}
{\it dedicated to the memory of Gian-Carlo Rota}
\end{center} 

\vspace{0.5in}

{\sc Abstract:} 
Two of the pillars of combinatorics
are the notion of choosing an arbitrary subset 
of a set with $n$ elements (which can be done in $2^n$ ways),
and the notion of choosing a $k$-element subset 
of a set with $n$ elements (which can be done in $n \choose k$ ways).
In this article I sketch the beginnings
of a theory that would import these notions
into the category of hedral sets in the sense of Morelli
and the category of polyhedral sets in the sense of Schanuel.
Both of these theories can be viewed as extensions
of the theory of finite sets and mappings between finite sets,
with the concept of cardinality being replaced by
the more general notion of Euler measure
(sometimes called combinatorial Euler characteristic).
I prove a ``functoriality'' theorem (Theorem 1)
for subset-selection
in the context of polyhedral sets,
which provides quasi-combinatorial interpretations
of assertions such as $2^{-1} = \frac12$
and ${1/2 \choose 2} = -\frac18$.
Furthermore, the operation of forming a power set
can be viewed as a special case of the operation
of forming the set of all mappings from one set to another;
I conclude the article
with a polyhedral analogue of the set of all mappings
between two finite sets,
and a restrictive but suggestive result (Theorem 2)
that offers a hint of what a general theory of
exponentiation in the polyhedral category
might look like.
(Other glimpses into the theory may be found in {\bf [Pro]}.)

{\large
\bigskip
\noindent
1. Introduction
\bigskip
\normalsize}

In the work of Hadwiger {\bf [Ha1] [Ha2]}, 
Lenz {\bf [Len]}, McMullen {\bf [McM]}, Rota {\bf[KR]}, 
Chen {\bf [Ch1] [Ch2]}, and others, 
the Euler characteristic of polytopes in $\R^n$ 
is extended to a valuation $\chi$
(the Euler measure)
on the algebra of sets generated by these polytopes.  
Included in this algebra of sets
(called ``hedral sets'' by Morelli {\bf[Mor]}
and ``polyhedral sets'' by Schanuel {\bf [Sch]})
are the finite subsets of $\R^n$,
for which the notion of Euler measure
and the notion of cardinality coincide.
This suggests the viewpoint 
(forcefully advocated by Schanuel)
that we view polyhedral sets and their $\Z$-valued Euler measures
as a generalization of 
finite sets and their $\N$-valued cardinalities.
Indeed, several of the key pillars 
on which combinatorics rests ---
concepts of union, Cartesian product and
selection of a subset of fixed size ---
extend to the polyhedral setting straightforwardly,
and this extension provides a natural interpretation
for assertions like ${-2 \choose 3} = -4$
(see section 2).

What else constitutes the foundation of combinatorics?
Many combinatorialists would point to
the notion of the set of all maps
(or injections, or surjections, or bijections)
from one finite set to another;
see for instance the ``twelve-fold way'' of Rota 
(so dubbed by Spencer) as described in Stanley's book {\bf [Sta]}.
As a special case,
the set of all maps from a finite set to the set $\{0,1\}$
corresponds naturally with the set of all subsets of the given set.

I will show that a polyhedral analogue of this foundation can be built,
using Morelli's notion of hedral maps
(piecewise-constant maps with finite range
whose level sets are polyhedral sets)
and Schanuel's more inclusive notion of polyhedral maps
(maps whose graphs are polyhedral sets).
Since in neither case
is the set of {\it all\/} such maps from $A$ to $B$
itself a polyhedral set
(except when $A$ and $B$ are finite or $A$ or $B$ is empty),
we cannot simply take Euler measure in the standard way;
but, at least in the case where the polyhedral set $A$
is a subset of $\R$,
we can divide this large set of maps into
finite-dimensional strata
by limiting the number of points of discontinuity/non-differentiability
and use this stratification to define a power series in $t$
whose value at $t=1$ ``should'' be the Euler measure
of the set of maps.
In many cases the series possesses
an analytic continuation
with a well-defined value at $t=1$,
which I call the {\it regularized Euler measure}.

For instance, let $A=(0,1)$ and $B=\{0,1\}$,
with $\chi(A)=-1$ and $\chi(B)=2$;
then we will see in section 5 that 
$B^A$, defined as the set of all polyhedral maps from $A$ to $B$,
has regularized Euler measure $\frac12$ ,
and the set of all two-element subsets of $B^A$ 
has regularized Euler measure $-\frac18$.
The preceding assertions correspond,
in some way not yet formalized,
to the algebraic assertions
that $\frac12 = 2^{-1} = \chi(B)^{\chi(A)}$
and $-\frac18 = {1/2 \choose 2}$
(under the interpretation of $x \choose k$
as a polynomial of degree $k$ in $x \in \R$).

A different theory of exponentiation of polyhedral sets
arises from focusing on power sets,
and making the restriction that one looks at
only the {\it finite\/} subsets of a given polyhedral set.
Since this theory is simpler
and has been developed further,
it is discussed first, in sections 3 and 4.
Then, in section 5, I treat
the set of maps from one polyhedral set to another,
under both Morelli's condition of hedrality
and Schanuel's weaker condition of polyhedrality.

I should warn the reader not to expect too much;
the two theorems proved in this article are really 
more like examples than like theorems ---
and what is more, I do not give a clear statement
of what sorts of definitions and theorems
these examples are special cases {\it of\/}.
It is my hope that at least some readers
will see that there are issues here that are
worthy of being pursued,
and that there is some ``spiritual'' satisfaction
in the prospect of a theory 
that might be a model of the rational numbers
in the same way that ordinary combinatorics
is a model of the non-negative integers.

I thank Scott Axelrod, John Baez, Beifang Chen,
Timothy Chow, Ezra Getzler, Greg Kuperberg, Michael Larsen, 
Ayelet Lindenstrauss, Haynes Miller, Lauren Rose, 
and Dylan Thurston for helpful conversations.
I also thank Joseph Kung for soliciting this article
and putting this special issue together,
and the anonymous referee for numerous helpful suggestions.

I acknowledge a special debt to Gian-Carlo Rota,
who persistently encouraged the investigations described here,
despite their inconclusive status both in 1996
(when I described my work to him)
and now (after years of neglect).
I also want to acknowledge the formative role
he played in my career as the author of the beautiful book
{\it Finite Operator Calculus\/},
which was part of what inspired my earliest explorations
of combinatorics.
His sagacity and warmth are missed by me and many others.

{\large
\bigskip
\noindent
2. Background
\bigskip
\normalsize}

Following Schanuel,
we define a {\it polyhedral set\/} in $\R^n$
as one that can be defined by a Boolean predicate
whose primitive terms are linear equations and inequalities
in $n$ variables.
Equivalently, the algebra of polyhedral sets
is the algebra generated by
the closed half-planes in $\R^n$
under intersection, union and complementation.
Note that polyhedral sets can be unbounded;
for instance,
$\{(x,y): |x|+|y| \leq 1\} \cup \{(x,0): x > 2\} \cup \{(5,5)\}$
is a polyhedral subset of $\R^2$.

There exists a unique function $\chi(\cdot)$ 
from the polyhedral algebra to the integers
such that 
\begin{equation}
\label{valuation}
\chi(A \cup B) = \chi(A) + \chi(B) - \chi(A \cap B)
\end{equation}
for all polyhedral sets $A,B$
and such that
$$
\chi(A) = (-1)^k
$$
if the polyhedral set $A$ is homeomorphic to $\R^k$
or equivalently is homeomorphic to the interior of a $k$-dimensional disk 
(see e.g.\ {\bf [Ch2]}).
For instance, the polyhedral set 
described at the end of the preceding paragraph 
has $\chi$-value $1-1+1=1$.
The function $\chi$, hereafter called simply Euler measure,
is invariant under homeomorphism,
and in the case of compact polyhedral sets,
it coincides with the alternating sum of the Betti numbers.
However, unlike the notion of Euler characteristic favored by topologists,
it is not a homotopy invariant
(e.g., $\chi$ of an open interval in $\R$ is $-1$
whereas $\chi$ of a closed interval is $+1$).
Property (\ref{valuation}) extends naturally
to an inclusion-exclusion principle:
$$\chi(A_1 \cup A_2 \cup \dots \cup A_m) = \ \ \ \ \ \ \ \ \ \ \ \ \ \ \ \ 
\ \ \ \ \ \ \ \ \ \ \ \ \ \ \ \ \ \ \ \ \ \ \ \ \ \ \ \ \ \ \ \ $$
$$\sum_{1 \leq i \leq m} \chi(A_i) -
\sum_{1 \leq i < j \leq m} \chi(A_i \cap A_j) + 
\dots  \pm \chi(A_1 \cap A_2 \cap \dots \cap A_m).$$

For polyhedral sets $A \subset \R^m$ and $B \subset \R^n$,
one obtains a Cartesian product $A \times B$ in $\R^{m+n}$
in the obvious way,
and it can easily be shown that $\chi(A \times B) = \chi(A) \chi(B)$.
Less trivially, given a polyhedral set $A \subset \R^n$
and an integer $k \geq 0$,
there is a way of defining a polyhedral set
${A \choose k} \subset (\R^n)^k$
that contains exactly one point in each free orbit
in the $k$-fold Cartesian product $A^k$
under the permutation action of the symmetric group $S_k$
(where the free orbits are those of cardinality $k!$,
i.e., the orbits of $n$-tuples of distinct real numbers);
this construction is well-defined up to polyhedral isomorphism (defined below)
and has the property that $\chi({A \choose k}) = {\chi(A) \choose k}$.
Since the operation $A \mapsto {A \choose k}$
(mapping Schanuel's category of polyhedral sets to itself)
seems to be new,
it is necessary to amplify on the preceding sentences
in the following paragraphs.

Following Schanuel, we define
a {\it polyhedral map\/} from a polyhedral set $A \subset \R^m$
to a polyhedral set $B \subset \R^n$
as a map $A \rightarrow B$ 
whose graph in $\R^{m+n}$ is itself a polyhedral set;
that is, a polyhedral map is a piecewise-affine map
that need not be continuous.
A {\it polyhedral isomorphism\/} is a bijective polyhedral map,
such as the map
$$
f(x) = \left\{ \begin{array}{ll}
	x-1 & \mbox{if $-1 < x < 0$} \\ 
	x+5 & \mbox{if $x = 0$} \\ 
	2x+1 & \mbox{if $0 < x < 1$}
	\end{array} \right.
$$
between the interval $(-1,1) \subset \R$
and the set $(-2,-1) \cup (1,3) \cup \{5\} \subset \R$.
It can be shown that two polyhedrally isomorphic polyhedral sets
have the same Euler measure.
On the other hand, following Morelli, we define
a {\it hedral map\/} from a polyhedral set $A \subset \R^m$
to a polyhedral set $B \subset \R^n$
to be a function $f:A \rightarrow B$
with finite range
whose level-sets are all polyhedral.

In both cases, we identify a function $f: \R^m \rightarrow \R^n$
with its graph in $\R^{m+n}$.

Note that if $A$ is finite,
the hedral and polyhedral maps from $A$ to $B$ 
are just the set-theoretic maps,
so that, in the case where both $A$ and $B$ are finite,
the set of hedral or polyhedral maps from $A$ to $B$
has cardinality $|B|^{|A|}$
(or, in other terms, Euler measure $\chi(B)^{\chi(A)}$).
If $A$ is infinite, then
most set-theoretic maps from $A$ to $B$ 
are neither hedral nor polyhedral.

Also note that when $B$ is finite,
the the set of hedral maps from $A$ to $B$
coincides with the set of polyhedral maps from $A$ to $B$;
when $B$ is infinite,
every hedral map is polyhedral
but not conversely.

We now define a {\it polyhedral ordering\/} 
on a polyhedral set $A \subset \R^n$
as a (total) linear ordering
whose graph in $A \times A$
is polyhedral;
an example of such an ordering
is the one induced by lexicographic ordering
of $n$-tuples in $\R^n$.
Henceforth, we restrict attention to lexicographic ordering.
We form ${A \choose k} \subset A^k$
by selecting, from each free orbit
in the $k$-fold Cartesian product $A^k$
under the permutation action of the symmetric group $S_k$,
the one $k$-tuple that is lexicographically least
in the induced ordering on $A^k$.
That is, $A \choose k$ consists of all $k$-tuples
$(v_1,v_2,\dots,v_k) \in A^k$
with $v_1 < v_2 < \dots < v_k$
in the lexicographic sense.
(Note that any polyhedral set
that contains exactly one point in each free orbit
in $A^k$ under the action of $S_k$
is canonically equivalent to the one defined above,
via a polyhedral isomorphism.)

Take for instance $A = (0,1) \cup (2,3) \subset \R$,
with $\chi(A) = -2$.
Then we find that $A \choose 3$ is a union of four sets:
$\{(x,y,z): 0<x<y<z<1\}$,
$\{(x,y,z): 0<x<y<1,\ 2<z<3\}$,
$\{(x,y,z): 0<x<1,\ 2<y<z<3\}$ and
$\{(x,y,z): 2<x<y<z<3\}$.
Each of these is homeomorphic to the interior of a 3-ball,
so $\chi({A \choose 3}) = (-1)^3 + (-1)^3 + (-1)^3 + (-1)^3 = -4$.
Note that this is also the value of $\chi(A) \choose 3$
where for each fixed $k \geq 0$ and all $x \in \R$
we define $x \choose k$ as $x(x-1)\cdots(x-k+1)/k!$.

To show that $\chi({A \choose k}) = {\chi(A) \choose k}$
for all polyhedral sets $A$ and all integers $k \geq 0$,
it suffices to show that the polyhedral set
$\{(v_1,v_2,\dots,v_k) \in A^k: v_i \neq v_j \ \mbox{for all $i \neq j$}\}$
has Euler measure $\chi(A) (\chi(A)-1) \cdots (\chi(A)-k+1)$,
since this larger polyhedral set is the disjoint union of
$k!$ polyhedrally-isomorphic copies of $A \choose k$.
One way to prove this is by iterated application of a Fubini theorem
(Theorem 3 in Morelli):
If a polyhedral set $A \subset B \times C$
has the property that every $C$-fiber has Euler measure $m$,
then $A$ has Euler measure $\chi(B) \cdot m$.
For example, in the case $k=2$,
it is clear that $\{(v_1,v_2) \in A^2: v_1 \neq v_2\}$
projects to its first factor with fibers of Euler measure
$\chi(A)-1$,
so by the Fubini principle its Euler measure must be
$\chi(A) (\chi(A) - 1)$;
the same approach works when applied iteratively
in the case of larger values of $k$.
This proof has the virtue of being entirely parallel
to the standard way of counting ordered $k$-element subsets of a finite set,
and it makes it clear where the factors $\chi(A)-i$ come from.
However, as an intimation of what is to come,
it is more helpful to prove the result by 
poset-theoretic inclusion-exclusion, as we now do.
(See Stanley's book {\bf [Sta]} for background on partially ordered sets.)

The {\it partition lattice\/} $\Pi_k$
is the set of unordered partitions of a $k$-element set,
ordered so that its top element is the 1-block partition $\hat{1}$
and its bottom element is the $k$-block partition $\hat{0}$.
It is easy to see that for any finite set $S$,
the product set $S^n = \{(s_1,s_2,\dots,s_n)\}$ 
splits up into $n$ disjoint subsets,
where the $k$th subset (with $k$ between 1 and $n$ inclusive)
consists of those $n$-tuples that contain 
exactly $k$ distinct elements of $S$
and therefore has cardinality $|S|^k$.
Ordinary inclusion-exclusion applied to the lattice $\Pi_k$
lets us turn this around and conclude that
the number of ordered $k$-element subsets of $S$ is 
$$\sum_{\pi \in \Pi_k} \mu(\hat{0},\pi) |S|^{N(\pi)},$$
where $\mu(\cdot,\cdot)$ is the M\"obius function of $\Pi_k$
and $N(\pi)$ is the number of blocks in the partition $\pi$.
(See exercise 44 of {\bf [Sta]}.)
Proceeding analogously in the case of polyhedral sets,
we see that for any polyhedral set $A$,
the product polyhedron $A^n = \{(v_1,v_2,\dots,v_n)\}$ 
splits up into $n$ disjoint polyhedral sets,
where the $k$th polyhedral set (with $k$ between 1 and $n$ inclusive)
consists of those $n$-tuples that contain 
exactly $k$ distinct elements of $A$
and therefore has Euler measure $\chi(A)^k$.
As before, inclusion-exclusion applied to the lattice $\Pi_k$
tells us that the Euler measure of $A \choose k$ is 
$$\sum_{\pi \in \Pi_k} \mu(\hat{0},\pi) \chi(A)^{N(\pi)}.$$
Since the two polynomials are the same
(aside from the fact that one has $\chi(A)$ where the other has $|S|$),
and since we know that the former is equal to $|S|(|S|-1)\cdots(|S|-k+1)$,
the desired result follows.

(The operations $\lambda^k: A \mapsto {A \choose k}$ are akin to
the $\lambda$-ring structures that Morelli defines
on the abelian group of hedral functions,
in Theorem 13 of his article and the Remark that follows it.)

As an aside,
it is interesting to note that although finding a polyhedral analogue
of $n \mapsto 2^n$ seem to require us to go beyond the polyhedral category
and into an extension involving infinite-dimensional objects,
finding a polyhedral analogue of the Fibonacci numbers does not.
Given a polyhedral set $P \subset \R$ of Euler measure $n$
(positive, negative, or zero),
the set of all finite subsets $S \subset P$ such that 
$\chi((P \setminus S) \cap (t,t'))$ is even
for all $t,t' \in S \cup \{-\infty,+\infty\}$
has Euler measure $(\phi^{n+1}-(-1/\phi)^{n+1})/\sqrt{5}$
\ (with $\phi=(1+\sqrt{5})/2$).

{\large
\bigskip
\noindent
3. The Small Power Set
\bigskip
\normalsize}

To motivate our definition
of the Euler series,
we begin by considering
the ``small power set'' discussed earlier,
that is, the 
set of all finite subsets of $A$,
where $A$ is a polyhedral set.
We denote the small power set of $A$ by $2^A$.
$2^A$ is not itself a polyhedral set,
but it can be realized as a formal union of disjoint polyhedral sets
in spaces of increasing dimensionality.
For instance, 
in the case where $A$ is an open (not necessarily bounded) interval,
we can view $2^A$ as consisting of a point (the empty subset of $A$),
an open interval (the set of singleton subsets of $A$),
an open triangle (the set of 2-element subsets of $A$),
an open 3-simplex (the set of 3-element subsets of $A$),
and so on.
These constituent pieces have Euler measure
$1,-1,1,-1,\dots$ respectively.
The terms do not yield a summable sequence,
but there is no doubt what 
the incorrigible formalist Euler would have done:
the series $1-1+1-1+\dots$ is after all a geometric series,
so it can be assigned a notional sum $1/(1-(-1))=1/2$.
More precisely, and more generally,
given a polyhedral set $A$
one can define a formal power series $E(t)$
in which the coefficient of $t^k$
is the Euler measure of the set of $k$-element subsets of $A$.
That is, $E(t)={m \choose 0}+{m \choose 1}t+{m \choose 2}t^2 + \dots$
where $m=\chi(A)$.
This is just the binomial expansion of $(1+t)^m$;
when the integer $m$ is non-negative, the series is a polynomial
whose value at $t=1$ is $2^m$,
while for $m<0$ the series converges in an open disk of radius 1
to a function that extends to the whole complex plane away from $t=-1$,
and whose value at $t=1$ is $2^m$.
Either way, the value that one gets for $\chi(2^A)$ is $2^{\chi(A)}$.
This ``functoriality'' is a sign that we are on the right track.
(Here and elsewhere in the article, functoriality is said to hold
when operations in the polyhedral category
carry over, under the mapping $\chi$,
to operations in the category of natural numbers,
as described by Schanuel {\bf [Sch]}.)

Some less trivial corroboration of the validity of our approach 
comes from looking at $B \choose 2$,
where $B$ is the power set of an open interval $A$.  
That is, we are looking at all the ways of choosing
two {\it distinct\/} finite subsets $S_1$, $S_2$ of the interval $A$.
Given $\sigma = \{S_1,S_2\} \in {B \choose 2}$, so that $S_1 \neq S_2$,
define the {\it support-set\/} (or {\it support\/})
of $\sigma$ as $S_1 \cup S_2$, and define the
{\it rank\/} of $\sigma$ as $|S_1 \cup S_2|$.  We define an
Euler series $E(t)$ in which the coefficient of $t^k$ is the
Euler measure of the rank-$k$ part of $B \choose 2$.
Given a $k$-element subset $S$ of $A$,
there are $3^{k-1}+3^{k-2}+\dots+3+1$ ways
to define $\sigma \in {B \choose 2}$ with support $S$.
(Namely: for each $i$ between $0$ and $k-1$ inclusive,
one can assign the first $i$ elements of $S$ to both $S_1$ and $S_2$,
the next element of $S$ to $S_1$ but not $S_2$,
and each of the next $k-i-1$ elements of $S$
to $S_1$ only, $S_2$ only, or both $S_1$ and $S_2$,
as one wishes.
All $3^{k-1}+3^{k-2}+\dots+3+1$ such pairs $S_1,S_2$ are distinct
from one another.
On the other hand, every pair of distinct finite subsets of $B$
whose union has cardinality $n$
can be written in this form in a unique way,
by labelling the two sets as $S_1$ and $S_2$
so that the smallest element in their symmetric difference
belongs to $S_1$.)
Thus the Euler measure of the rank-$k$ part of $2^A$
is $\chi({A \choose k}) \cdot (3^{k-1}+3^{k-2}+\dots+3+1)$.
The resulting Euler series $-t+(3+1)t^2-(9+3+1)t^3+\dots$
analytically extends to the rational function $-t/(1+t)(1+3t)$,
whose value at $t=1$ is $-1/8 = {1/2 \choose 2}$.

Let us agree to limit ourselves to taking power sets
of subsets of $\R^2=\R$.
(Note that the Euler measure of a polyhedral subset of $\R$
can be any integer, 
so not much generality has been lost through this restriction,
in terms of the sort of numbers that will arise
as Euler measures of derived sets.)
We are going to consider sets of the form
$${B \choose k_1},$$
$${{B \choose k_1} \choose k_2},$$
$${{{B \choose k_1} \choose k_2} \choose k_3},$$
etc., where $A$ is a polyhedral subset of $\R$
and $B=2^A$;
we will define them iteratively.
We will call these sets {\it gizmos\/}
(a more respectable term might be ``selection schemes'').
At each step, we will introduce a linear ordering on the gizmo,
along with a support-set operation;
together these allow us to proceed to the next step.
Our notation for such gizmos
will be $G(2^A;k_1)$,
$G(2^A;k_1,k_2)$,
$G(2^A;k_1,k_2,k_3)$, etc.

We start with $2^A$ itself.
Given subsets $S,T$ of $A$
with $|S| = |T|$,
we say that $S$ is less than $T$
if the smallest element in their symmetric difference
belongs to $S$.
(Equivalently: Writing $S = \{s_1,s_2,\dots,s_k\}$
and $T = \{t_1,t_2,\dots,t_k\}$
with $s_1<s_2<\dots<s_k$
and $t_1<t_2<\dots<t_k$,
$S$ is less than $T$ if $s_i < t_i$
for the first location $i$ at which the two sequences differ.)
The support-set operation from $2^A$ to itself
is just the identity operation.
Given a gizmo $G$ with a linear ordering and a support-set operation,
and given an integer $k \geq 0$,
we define the gizmo $G \choose k$ as the set of all sequences
$(g_1,g_2,\dots,g_k) \in G^k$ with $g_1 < g_2 < \dots < g_k$
(with respect to the linear ordering on $G$).
We order $G \choose k$ lexicographically
(so that $(g_1,g_2,\dots,g_k)$ is less than $(g'_1,g'_2,\dots,g'_k)$
if and only if $g_i < g'_i$ at the first location $i$ at which
the two sequences differ)
and we define the support of $(g_1,\dots,g_k)$ as
the union of the supports of $g_1$ through $g_k$.

For any gizmo $G$, the number of elements of $G$
whose support is a particular $k$-element subset of $A$
does not depend on the particular subset of $A$,
but is some number $n_k$ that depends only on $k$ (and $G$).
We therefore define the Euler series of $G$
as
$$\sum_{k=0}^\infty \chi({\scriptstyle {A \choose k}}) n_k t^k.$$

\vspace{0.2in}

{\large
\bigskip
\noindent
4. Continued Binomial Coefficients
\bigskip
\normalsize}

{\sc Theorem 1:} The Euler series for 
$G(2^A;k_1,k_2,\dots,k_r)$ converges in a neighborhood of $t=0$
to a rational function of $t$
whose value at $t=1$ is the iterated binomial coefficient
$C(2^{\chi(A)};k_1,k_2,\dots,k_r)$
(defined by the initial condition $C(n)=n$
and the recurrence relation $C(n;k_1,k_2,\dots,k_i)=
{C(n;k_1,k_2,\dots,k_{i-1}) \choose k_i}$).

{\sc Remark:}
While this result is not as general as one might like,
it does back up my contention
that my definition of the power set,
or something like it,
can interact in a functorially appropriate way
with other more basic polyhedral operations,
once the category of polyhedral sets and polyhedral maps
has been extended in a suitable way.

{\sc Proof:}
As in the argument given at the end of Section 2,
the proof of the main theorem uses inclusion-exclusion.
We illustrate the proof for the concrete case
$G(2^{\chi(A)};2,2)$
where $A$ varies over the polyhedral subsets of $\R$,
since it contains all the key ideas.
An element of this gizmo can be described
as a 4-tuple $(S_1,S_2,T_1,T_2)$ of subsets of $A$
satisfying the following three conditions:
(1) $S_1 < S_2$;
(2) $T_1 < T_2$;
(3) $\{S_1,S_2\} < \{T_1,T_2\}$.
Note that condition (3) is equivalent to the condition
that either (i) $S_1 < T_1$ or else (ii) $S_1 = T_1$ and $S_2 < T_2$.
We know that, in the case where $A$ is a finite set of points,
the gizmo is itself just a finite set of points,
and that its Euler series is a polynomial in $t$
whose value at $t=1$ is just 
$C(2^{\chi(A)};2,2)$.
Our goal is to show that this same formula 
applies even when $A$ is not a finite set of points,
essentially by showing that the two cases are governed by
the same inclusion-exclusion calculation.

Before we do the comparison, it is helpful
to switch to a re-statement of the problem
more amenable to inclusion-exclusion analysis.
We will the consider the ``symmetrized gizmo''
consisting of all 4-tuples obtained from the elements of
$G(2^{\chi(a)};2,2)$ by the action of the 8-element group that
allows one to switch $S_1$ with $S_2$, $T_1$ with $T_2$,
or the pair $S_1,S_2$ with the pair $T_1,T_2$.
Note that the elements of the symmetrized gizmo
are exactly the 4-tuples $(S_1,S_2,T_1,T_2)$ of subsets of $A$
for which {\it none\/} of the following conditions holds:
(a) $S_1 = S_2$;
(b) $T_1 = T_2$;
(c) $S_1 = T_1$ and $S_2 = T_2$;
(d) $S_1 = T_2$ and $S_2 = T_1$.
Regardless of the nature of $A$ (finite vs.\ infinite),
it is clear that the Euler series of the new gizmo
is 8 times the Euler series of the old,
so it is enough to prove that the Euler series of the symmetrized gizmo
is a rational function of $t$ whose value at $t=1$
is given by a polynomial in $2^{\chi(A)}$.
Once we have shown that the value at $t=1$
depends only on $\chi(A)$, and not on $A$,
our observations for the case where $A$ is finite will
complete the proof.

To prove the required properties of the Euler series
for symmetrized gizmos,
the principle of inclusion-exclusion guarantees that
it is enough to show that the claim holds
when we look at the set of 4-tuples
in which some specified (possibly empty) subset
of the forbidden conditions (a),(b),(c),(d)
does in fact hold.
But the effect of imposing any or all of these four conditions
is merely to require that one is not choosing four subsets of $A$
but three, two, or one.
Specifically, if we let $E_{a}$
denote the Euler series for the set of 4-tuples of (finite) 
subsets of $A$ for which condition (a) holds, etc.,
and we let $E(k)$ denote the Euler series for
the set of $k$-tuples of subsets of $S$,
then the Euler series for $G(2^A;2,2)$
is equal to
$[E] - [E_a + E_b + E_c + E_d] 
+ [E_{ab} + E_{ac} + E_{ad} + E_{bc} + E_{bd} + E_{cd}]
- [E_{abc} + E_{abd} + E_{acd} + E_{bcd}] + [E_{abcd}]
= [E_4] - [E_3 + E_3 + E_2 + E_2] + [E_2 + E_1 + E_1 + E_1 + E_1 + E_1]
- [E_1 + E_1 + E_1 + E_1] + [E_1]
= E_4 - 2E_3 - E_2 + 2E_1$.
(Note that if we replace $E_k$ by $x^k$,
the preceding expression reduces to
$$8 {{x \choose 2} \choose 2},$$
as it ought to.)
Here we see 
that the M\"obius coefficients 
that enter into the inclusion-exclusion formula
are independent of the set $A$,
but depend only on $k_1 = 2$, $k_2 = 2$.
This independence is true more generally.
Hence the entire problem reduces to the case of $(2^A)^j$,
with $j$ varying;
that is, the gizmo consisting of all ways of choosing
a $j$-tuple of (not necessarily distinct) subsets 
$S_1,S_2,\dots,S_j$ of $A$.

We can now conclude the proof.
When $A$ is a finite set of points,
the Euler series for $(2^A)^j$ is
$(1+(2^j-1)t)^{\chi(A)}$,
since for each of the $\chi(A)$ points of $A$
we can either choose to include it in {\it none\/} of the sets
$S_1,S_2,\dots,S_j$
or else choose to include it in {\it some\/} of the $j$ sets
(which we can do in any of $2^j-1$ ways).
When $A$ is an open interval,
the Euler series for $(2^A)^j$ is
$(1+(2^j-1)t)^{-1}=1-(2^j-1)t+(2^j-1)^2t^2-\dots$,
since for each subset of $A$ of size $i$,
the number of ways of assigning each element of the subset
to at least one (and possibly all) of the $j$ sets $S_1,S_2,\dots,S_j$
is $(2^j-1)^i$.
In general,
the Euler series for $(2^A)^j$ is
$(1+(2^j-1)t)^{\chi(A)}$,
which when evaluated at $t=1$
gives $(2^j)^{\chi(A)} = (2^{\chi(A)})^j$.
Note that this is a polynomial in $2^{\chi(A)}$, as claimed.
Hence all the gizmos in question
(the original gizmo $G(2^A;k_1,k_2,\dots,k_r)$ 
and its symmetrized version)
have Euler measure given by polynomials in $\chi(A)$,
which we can identify by looking at the case where $A$
is a finite set.  $\Box$

\medskip

Theorem 1 in its current form is a bit eccentric in its focus.
Ideally one would like to replace it by a result
that says that the operation selecting a subset of size $k$
may be freely combined with the operations of
disjoint union and Cartesian product
(and perhaps other natural polyhedral operations as well)
in such a fashion that all the resulting polyhedral sets
satisfy formulas like
$\chi(A \choose k) = {\chi(A) \choose k}$,
$\chi(A \cup B) = \chi(A) + \chi(B) - \chi(A \cap B)$,
and $\chi(A \times B) = \chi(A) \chi(B)$.

{\large
\bigskip
\noindent
5. Exponentiation
\bigskip
\normalsize}

If $A,B$ are polyhedral sets
($A \subset \R^m$, $B \subset \R^n$),
the set of {\it affine\/} maps from $A$ to $B$
has a well-defined Euler measure
that can be associated with it in a number of ways.
One simple way of defining it
is to let $v_1,v_2,\dots,v_k$ be points in $\R^m$
whose affine span contains $A$;
then every map from $A$ to $B$
is determined by its values on $v_1,v_2,\dots,v_k$,
so that the set of affine maps from $A$ to $B$
becomes associated with a subset of $B^k$
(namely, the set of all $k$-tuples $(f(v_1),f(v_2),\dots,f(v_k))$
arising from affine maps $f: A \rightarrow B$).
One can show that this set in $B^k$ is itself polyhedral,
and that its Euler measure does not depend
on the spanning set $v_1,v_2,\dots,v_k$ that was chosen.

For example, in the case $A=(0,1) \subset \R$
it is natural to take $v_1=0$ and $v_2=1$.
Assume that $B$ is compact.
Then the set of affine maps from $A$ to $B$
corresponds to the set of pairs $(w,w') \in \R^n \times \R^n$
such that $tw+(1-t)w' \in B$ for all $0 < t < 1$.
It is not hard to show that this set is polyhedral
and that for each $w \in \R^n$,
the set $\{w' \in \R^n: \mbox{$tw+(1-t)w' \in B$ for all $0<t<1$}\}$
has Euler measure 1 or 0
according to whether or not $w \in B$
(this is where compactness gets used,
via the fact that a starlike compact polyhedral set
has Euler measure 1).
Hence, by Fubini's Theorem,
the set of affine maps from $(0,1)$ to $B$
has Euler measure $\chi(B)$.

Let us now take $A$ to be a polyhedral subset of $\R$
with no isolated points.
(It is clear that the case of general polyhedral sets $A$ can be handled
by combining the case in which $A$ is finite
and the case in which $A$ has no isolated points.)
Then any polyhedral map $f:A \rightarrow B$ is locally affine
at all but a finite number of points.
We call a point at which $f$ fails to be locally affine
(whether by virtue of discontinuity or by virtue of non-differentiability)
a {\it breakpoint\/} of $f$,
and we call the number of breakpoints
the {\it rank\/} of $f$.
(Note that this notion of rank, specialized
to the case $B=\{0,1\}$,
coincides with the notion introduced in section 3:
since we are assuming that $A$ has no isolated points,
the indicator function of a finite subset of $A$ with cardinality $n$
has exactly $n$ points of discontinuity,
namely, the elements of $A$ themselves.)

Let $S$ be a finite subset of $A$.
We wish to count (in the sense of Euler-measure)
the polyhedral maps from $A$ to $B$
whose breakpoints are precisely the points of $S$.
It is not hard to see that the set of maps from $A$ to $B$
that are locally affine on $A \setminus S$
is amenable to the same analysis as we applied
in the first paragraph of this section;
that is, it has a well-defined Euler measure.
If a set is locally affine on $A \setminus S$,
its breakpoints form a subset of $S$,
and vice versa.
Therefore the method of inclusion-exclusion
gives us a way to calculate the Euler measure
for the set of polyhedral maps from $A$ to $B$
with breakpoint-set $S$.
It is easy to prove that this Euler measure has the same value
for all $S$ of fixed size $k$.
Call this value $n_k$.
Then one can define the Euler measure
of the set of polyhedral maps with exactly $k$ breakpoints
as $\chi({A \choose k}) \cdot n_k
= {\chi(A) \choose k} n_k$,
and define $E(t) = \sum_{k=0}^\infty {\chi(A) \choose k} n_k t^k$.

For example, suppose $A=(0,1)$, $B=\{0,1\}$.
A polyhedral map from $A$ to $B$ will have some number of
breakpoints, say at $0<x_1<x_2<\dots<x_k$ ($k \geq 0$);
in between these breakpoints, the map will be constant,
since $B$ is discrete.
I claim that the number of polyhedral maps from $A$ to $B$
with some fixed set of $k$ breakpoints is $2 \cdot {3^k}$.
For, starting at the left, the map can take either the value
0 or the value 1 on the interval $(0,x_1)$.  Thereafter,
for each $1 \leq i \leq k$, we can choose either value (0 or 1)
for the function at $x_i$ and either value for the function
on $(x_i,x_{i+1})$, as long as we do not have {\it both\/}
values equal to the value on $(x_{i-1},x_i)$ (for then
there would not be a breakpoint at $x_i$); there are
always 3 ways to implement this.  Hence the Euler series
is $2-6t+18t^2-54t^3+\dots=2/(1+3t)$, whose value at $t=1$
is $2/4=2^{-1}=\chi(B)^{\chi(A)}$.

Continuing this example in the spirit of Theorem 1,
we can look at the set of pairs of distinct maps 
from $(0,1)$ to $\{0,1\}$,
where the rank of such a pair is the cardinality of
the union of the sets of breakpoints of the two maps.
Then one can check that the set of maps
with $k$ specified breakpoints
has Euler measure
$2(2^4-1)^k-(2^2-1)^k = 2 \cdot 15^k - 3^k$,
so that $E(t)=\frac{2}{1+15t}-\frac{1}{1+3t}$
and $E(1)=\frac{2}{16}-\frac{1}{4}=-\frac{1}{8}={1/2 \choose 2}$
as expected.

Let us now focus on the case
where $A=(0,1)$ and $B$ is a bounded polyhedral set.
There are two sorts of exponentiation we will consider:
the ``Morelli exponential''
(the set of all hedral maps from $A$ to $B$)
and the ``Schanuel exponential''
(the set of all polyhedral maps from $A$ to $B$);
for definitions, see section 1.
The two cases are similar,
for in both cases,
the set of affine maps from any sub-interval $(x_i,x_{i+1})$ to $B$
has Euler measure $\chi(B)$.
(In the Morelli case, this is trivial,
since one is merely looking at constant maps;
in the Schanuel case, one appeals to
the analysis given in the second paragraph of this section.)
This allows us to prove:

{\sc Theorem 2:} Let $B$ be any bounded polyhedral set, and
let ${\cal F}$ be the set of hedral maps from $(0,1)$ to $B$
or the set of polyhedral maps from $(0,1)$ to $B$.
In either case, the regularized Euler measure of ${\cal F}$
is 0 if $\chi(B)=0$ and $1/\chi(B)=\chi(B)^{-1}$ otherwise.

{\sc Proof:}
For every set $S$ of $k$ points in $(0,1)$,
the set of hedral or polyhedral maps from $(0,1)$ to $B$
whose breakpoints lie in $S$
has Euler measure $\chi(B)^{2k+1}$,
with $k$ factors of $\chi(B)$ coming from the points in $S$
and $k+1$ factors coming from the intervals in between
(by virtue of the fact mentioned in the paragraph
preceding the statement of Theorem 2).
Inclusion-exclusion on the lattice of subsets of $S$
implies that the set of maps from $(0,1)$ to $B$
whose breakpoint-set is $S$
has Euler measure $\chi(B) (\chi(B)^2-1)^k$.
Hence $E(t)=\chi(B) / (1+(\chi(B)^2-1)t)$.
If $\chi(B)=0$, all coefficients vanish;
if $\chi(B) \neq 0$, $E(1) = \chi(B) / (\chi(B)^2) = 1/\chi(B)$.
$\Box$

{\large
\bigskip
\noindent
6. Final thoughts 
\bigskip
\normalsize}

It is worth pointing out here that full functoriality
is too much to expect; for instance, if 
$A$ has generalized Euler measure $1/2$
and $B$ has Euler measure $-1$, then
it is hard to believe that one could come up with a theory that
would assign to $C=B^A$ a generalized Euler measure
$(-1)^{\frac12}=\pm i$, and still harder to believe such a theory
would know what to make of $C^C$ (given that $i^i$ takes
on infinitely many values).
Yet numerous examples demonstrate
that as long as one does not push things too far,
functoriality does hold.

I do not claim that my proposed definition
of regularized Euler measure
is the ``right'' one;
if the reader finds its formulation artificial,
I must to some extent agree.
However, I think that the way the theory-in-progress hangs together
indicates that there is probably a better-motivated definition
that agrees with mine for the class of examples I consider.
A homological formulation would be especially useful,
and I hope some readers will pursue this line of investigation.

\vspace{0.4in}

{\large
\bigskip
\noindent
References
\medskip
\normalsize}

\noindent
{\bf [Ch1]} \ \ B.\ Chen, The Gram-Sommerville and Gauss-Bonnet theorems 
and combinatorial geometric measures for noncompact polyhedra,
{\it Advances in Mathematics} {\bf 91} (1992), 269-291.

\medskip

\noindent
{\bf [Ch2]} \ \ B.\ Chen, On the Euler measure of finite unions of 
convex sets,
{\it Discrete and Computational Geometry} {\bf 10} (1993), 79-93.

\medskip

\noindent
{\bf [Ha1]} \ \ H.\ Hadwiger, Eulers Characteristik und kombinatorische
Geometrie,
{\it J.\ reine angew.\ Math.\/} {\bf 194} (1955), 101--110.

\medskip

\noindent
{\bf [Ha2]} \ \ H.\ Hadwiger, Erweiterter Polyedersatz und 
Euler-Shephardsche Additionstheoreme, 
{\it Abhandlungen Mat.\ Seminar Hamburg\/} {\bf 39}
(1973), 120--129.

\medskip

\noindent
{\bf [KR]} \ \ D.\ Klain and G.-C.\ Rota,
Introduction to geometric probability.
Cambridge University Press, 1997.

\medskip

\noindent
{\bf [Len]} \ \ H.\ Lenz, Mengenalgebra und Eulersche Charakteristik,
{\it Abh.\ Math.\ Sem.\ Univ.\ Hamburg} {\bf 34} (1970), 135--147.

\medskip

\noindent
{\bf [McM]} \ \ P.\  McMullen, The polytope algebra, 
{\it Advances in Mathematics} {\bf 78} (1987), 76--130.

\medskip

\noindent
{\bf [McS]} \ \ P.\ McMullen and R.\ Schneider, Valuations on convex 
bodies, {\it Convexity and its Applications,}
P.\ Gruder and J.\ Wills, eds., Birkh\"auser, Basel 1983, pp.\ 160--247.

\medskip

\noindent
{\bf [Mor]} \ \ R.\ Morelli, A theory of polyhedra, {\it Advances in 
Mathematics}
{\bf 97} (1993), 1--73.

\medskip

\noindent
{\bf [Pro]} \ \ J.\ Propp, Euler measure as generalized cardinality;
arXiv: {\tt math.CO/}
{\tt 0203289}. 

\medskip

\noindent
{\bf [Sch]} \ \ S.\ Schanuel, Negative sets have Euler characteristic
and dimension, in: Lecture Note In Mathematics \#1488, Springer, 1991.

\medskip

\noindent
{\bf [Sta]} \ \ R.\ Stanley, Enumerative combinatorics I,
Wadsworth, 1986.

\end{document}